\documentclass[letterpaper, 10 pt, conference]{ieeeconf}  

\IEEEoverridecommandlockouts                              
\usepackage{amsmath} 
\usepackage{amssymb}  
\newcommand\R{\mathbb{R}}
\newcommand\C{\mathbb{C}}
\newcommand\N{\mathbb{N}}

\renewcommand\cal{\mathcal}

\newtheorem{theorem}{Theorem}
\newtheorem{remark}{Remark}
\title{\LARGE \bf
Control through operators for quantum chemistry
}

\author{Philippe Laurent$^{1}$, Herschel Rabitz$^{2}$, Julien
  Salomon$^{3}$ and Gabriel Turinici$^{3}$
\thanks{$^{1}$P. Laurent is with the IRCCyN,
Ecole des Mines de NANTES, 4 rue Alfred Kastler, 44300 Nantes, FRANCE
{\tt\small philippe.laurent at mines-nantes.fr}}
\thanks{$^{2}$H. Rabitz is with Department of Chemistry, Princeton
  University, Princeton, NJ 08544, USA
        {\tt\small hrabitz at princeton.edu}}%
\thanks{$^{3}$J. Salomon and G. Turinici are with the CEREMADE,
  Universit{\'e} Paris-Dauphine, Pl. du Mal. de Lattre de Tassigny,
  75016 Paris, FRANCE  {\tt\small julien.salomon at dauphine.fr, garbiel.turinici at dauphine.fr}}%
}

\begin{document}

\maketitle
\thispagestyle{empty}
\pagestyle{empty}

\begin{abstract}
We consider the problem of operator identification in quantum
control. The free Hamiltonian and the dipole moment are searched such
that a given target state is reached at a given time. A local
existence result is obtained. As a by-product, our works reveals
necessary conditions on the laser field to make the identification
feasible. In the last part of this work, some Newton algorithms are
proposed together with a continuation method to compute effectively
these operators. 

\end{abstract}

\section{Introduction}

In the last decades, quantum control has known significant improvements
both at theoretical and practical levels (cf.\cite{hbref1,hbref4,hbref9,hbref10} and references therein). Results
have been obtained on existence of controls~\cite{B,BL,controllability1,albertini,altafini} or efficient ways to
compute and carry out laser fields that
achieve some goals concerning the state of quantum
systems~\cite{BW,MT,palao,moi,MST2,MST}. On the other hand, the design of relevant laser
fields plays also a major role when the goal is to identify some
properties of the quantum system to be controlled. In this way, some methods have
been designed to identify finite dimensional
systems characteristics~\cite{LB}, or to compute discriminant laser
fields~\cite{cdc2009}.

Note that operator identification in relation to the Schr\"odinger
equation has already been studied in the literature. As an example, we refer
to~\cite{baudouin} for a theoretical result, where no laser
interaction is considered.

In this paper, we focus on the case where only one fixed laser is used to
identify in finite given time the free Hamiltonian and the dipole moment. From the
theoretical point of view, we obtain a local existence result:
we prove that the inversion is always possible in the neighborhood of some
particular states. As a by-product, we emphasize some features
of the laser fields that enables the identification. 

Following the local approach we use to obtain this result, we present
in a second part, a time discretized setting and fixed-point methods 
to solve numerically our problem. In particular, a Newton method is
proposed together with a continuation method that allows us to solve
problems where the local assumption does not hold.

This paper is organized as follows: the mathematical formulation of
our problem is given in Section~\ref{sec:setting} and a local
controllability result is presented in Section~\ref{sec:locrts}. In
Section~\ref{sec:algo}, we present the algorithms to solve numerically
the identification problem. We conclude with some tests in Section~\ref{sec:test}.

Let us finally introduce some notations concerning particular matrix
sets that will be used throughout the paper. Given
$N_d\in \N$, we denote by 
$\C^{N_d,N_d}$ and $\R^{N_d,N_d}$ the sets
of matrices of size $N_d\times N_d$ with complex and real coefficients respectively. Then,
define
\begin{eqnarray*}
\cal U\!\! &\!=\!&\!\!\left\{ M\in \C^{N_d,N_d}, \ M^*M=MM^*=Id \right\},\\
\cal S\!\! &\!=\!&\!\!\left\{ M\in \C^{N_d,N_d}, \ M^ *= M \right\},\\
\cal S_\R\!\! &\!=\!&\!\!\left\{ M\in \R^{N_d,N_d}, \ M^ *= M \right\},\\
\cal S^0_\R\!\! &\!=\!&\!\!\cal S_\R \cap \left\{ M\in
\R^{N_d,N_d},\ M_{k,k}=0,\  k=1,\ldots, N_d \right\},
\end{eqnarray*}
where $M^*$ denotes the adjoint matrix associated to $M$ and $Id$ is
the identity matrix of $ \C^{N_d,N_d}$. Here, for the sake of
simplicity, we have omitted the dependence of these sets with respect to
$N_d$. In what follows, we denote by $\Re z$ and $\Im z$ denote respectively the real and imaginary
parts of a complex number $z$. Given a matrix $M$, we denote by $M^T$ its transposed.

\section{Setting of the problem}\label{sec:setting}

 Fix $T>0$, and consider a system $U(t)\in \cal U$ whose
dynamics over $[0,T]$ is ruled by the Schr\"odinger equation:
\begin{eqnarray}
i\dot U(t)&=&[H_0+\varepsilon(t)\mu]U(t),\label{eq:schrod1}\\
U(0)&=&U_{init}\label{eq:schrod2},
\end{eqnarray}
where $H_0\in\cal S_\R$ is the matrix of the internal Hamiltonian, $\varepsilon(t)\in L^2(0,T;\R)$ a laser field,
$\mu\in \cal S_\R$ the matrix associated with the dipole moment. For relevant
applications, the matrices $H_0$ and $\mu$ are not supposed to
commute. The initial state $U_{init}$ is fixed.
In this equation, $\varepsilon$ is given and the pair $(H_0,\mu)\in\cal S_\R\times \cal S^0_\R$ is
searched such that at time $t=T$, the state reaches a given 
target state  $U_{target}$, i.e.,
\begin{equation}\label{eq:target} 
U(T)=U_{target}.
\end{equation}
In other words, given the mapping 
$$\begin{array}{cccl}
\varphi:& \cal S_\R\times \cal S^0_\R&\rightarrow& \cal U\\
&(H_0,\mu)&\mapsto& U(T),
\end{array} $$
the main question that will be investigated in this paper is the surjectivity of $\varphi$.

In our work, the internal Hamiltonian $H_0$ is searched as real Hermitian (i.e. symmetric) matrix. 
This is a particular situation as in general it is only supposed to be complex Hermitian and not real.
Nevertheless, for the applications we have in mind this restriction is very natural since the Hamiltonian is a sum of
a kinetic operator and a potential, both real. 
For the same reasons, we
suppose that the dipole moment $\mu$ is real (Hermitian thus symmetric) but
we assume moreover that the diagonal elements are null. 
This additional assumption is motivated both by invariance properties (the diagonal of $H_0$ as matrix commutes with the diagonal of $\mu$ as matrix)
but also by the desire to identify an unique pair $(H_0,\mu)$ since in this way the number of unknowns (dimension of $ \cal S_\R$ plus that of $\cal S^0_\R$) equals the  
number of equations (the dimension of $\cal U$). 

Note that one can easily prove the following conservation property:
$$\forall t\in[0,T],\ \|U(t)\|_{\cal U}=\|U_{init}\|_{\cal U},$$
where we have denoted by $\|\cdot \|_{\cal U}$ the norm associated to
the scalar product $$(A,B)\in {\cal U}\times{\cal U}\mapsto tr(A^* B).$$
This problem is related to inverse problems in quantum
control~\cite{cdc2009}, but unlike previous works, we do not
aim here at designing relevant laser fields to identify  the pair
$(H_0,\mu)$ but rather to investigate the properties of the fields
$\varepsilon(t)$ that make Equation~\eqref{eq:target} invertible and algorithms to compute
numerically the corresponding solution operators $H_0$ and $\mu$.

\section{Local controllability result}\label{sec:locrts}

In this section, we present some theoretical results about the local
inversion of Equation~\eqref{eq:target}. More precisely, we make use of the calculus of
variations to obtain a local inversion theorem.

Given a pair $(H_0,\mu)$, we first introduce the tangent space
$\cal A_{H_0,\mu}$, which is the space of matrices defined by:
$$\cal A_{H_0,\mu}=\left\{M \in \C^{N_d,N_d}, \ M^*U(T)+U(T)^*M=0 \right\}.$$
We then consider the differential operator of $\varphi$ defined by:
\begin{eqnarray*}
d\varphi (H_0,\mu):& \cal S_\R\times \cal S_\R\rightarrow \cal A_{H_0,\mu}\\
&(\delta H_0,\delta \mu)\mapsto \delta U(T),
\end{eqnarray*} 
where $\delta U(T)$ is solution at time $t=T$ of the linearized
Schr\"odinger equation:
$$i\dot{ \delta U}(t)=[H_0+\varepsilon(t)\mu]\delta U(t)+[\delta
  H_0+\varepsilon(t)\delta \mu]U(t),$$
and $U(t)$ follows equation~\eqref{eq:schrod1}. 

We will prove that $\varphi$ is an onto mapping using the fact
that $d\varphi$ also satisfies this property. This strategy is
motivated by the following known result:
\begin{theorem}
Supposed that $d\varphi (H_0,\mu)$ is an onto mapping, i.e.
$$\forall V\in \cal A_{H_0,\mu},\exists (\delta H_0,\delta\mu),
 d\varphi (H_0,\mu)(\delta H_0,\delta\mu)=V.$$
Then $\varphi$ is locally onto in a neighborhood  of $(H_0,\mu)$.
\end{theorem}
We shall prove that $d\varphi$ is an onto mapping on the neighborhood
of all states of the form $U_0:=\varphi (H_0,0)\in \cal U$. To do this,
we compute explicitly an inverse mapping.
\begin{theorem}\label{Th:ontoloc}
Given $H_0\in \cal S_\R$, define $V_0$ as the matrix that diagonalizes $U_0:=\varphi (H_0,0)$ in the
following way:
$$ U_0(t)=V^*_0e^{i\Lambda (t-\frac T2)}V_0,$$
with  $\Lambda$ the diagonal matrix with coefficients
$\lambda_a\in\R,\ a\in\N_d,\ 1\leq a\leq N_d$. Suppose that for $a\neq b,
1\leq a\leq N_d, 1\leq b\leq N_d$,
\begin{eqnarray}
\lambda_a&\neq&\lambda_b\label{cond1}\\
\hat{\varepsilon}^i_{a,b}&:=&\Im \left(\int_0^T
\varepsilon(t)e^{i\delta\lambda_{a,b}(t-\frac T2)}dt\right)\neq 0.\label{cond2}
\end{eqnarray} 
Then $d\varphi (H_0,0)$ is an onto mapping and
its inverse is given by:
$$\psi: V' \in\cal A_{H_0,\mu} \mapsto (\delta H_0,\delta \mu).$$
The matrices $\delta H_0$ and $\delta \mu$ are given by:
$$\delta H_0:=V_0^*\delta \tilde H_0 V_0,\ \delta \mu:=V_0^*\delta \tilde \mu V_0,$$
where the coefficients $h_{a,b}$
and $m_{a,b}$ of the matrices $\delta\tilde H_0$
and $\delta\tilde \mu$ are given by:

\begin{equation}
\left\{
\begin{array}{ccll}
m_{a,b}&=&\dfrac{\Im v_{a,b}}{\hat{\varepsilon}^i_{a,b}}&\\
h_{a,b}&=&\dfrac{\Re
  v_{a,b}-\frac{\hat{\varepsilon}^r_{a,b}}{\hat{\varepsilon}^i_{a,b}}\Im
  v_{a,b}}{\sin(\delta\lambda_{a,b}\dfrac T2)}\delta\lambda_{a,b} & {\textrm  if } a\neq b \\
m_{a,a}&=&0,&\\ 
h_{a,a}&=&\dfrac 2T v_{a,a} & {\textrm  if } a=b. 
\end{array}
\right.
\end{equation}
Here $v_{a,b}, \ a,b\in\N_d,\ 1\leq a\leq N_d$ are the coefficients of  $iV^*_0
U_0(T)^*V'V_0$ and $\hat{\varepsilon}^r_{a,b}:=\Re \left(\int_0^T
\varepsilon(t)e^{i\delta\lambda_{a,b}(t-\frac T2)}dt\right)$.
\end{theorem}

\begin{proof} We fix $V'\in \cal A_{H_0,\mu}$ and solve
\begin{equation}\label{eq:dphisurj}
d\varphi (H_0,0)(\delta H_0,\delta\mu)=V'.
\end{equation}
\\
First, one can show
the identities:
\begin{eqnarray}
 \varphi(H_0,\mu)^*d\varphi(H_0,0)(\delta H_0,\delta
 \mu)=U_0(T)^*\delta U_0(T)\nonumber\\
=-i\int_0^T  U_0(t)^*(\delta H_0+\varepsilon(t)\delta\mu)U_0(t)dt,\label{eq:idrem}
\end{eqnarray}
where the variation $\delta U_0$ is defined by the evolution equation:
\begin{equation}\label{eq:diffschrod}
i\dot{\delta U_0}(t)=[ H_0+\varepsilon(t)\mu]\delta U_0(t)+[ \delta
  H_0+\varepsilon(t)\delta \mu]U_0(t).
\end{equation}
Note that such an identity holds also when $\mu \neq 0$.
Since $U_0(T)^*$ is invertible, showing
that~\eqref{eq:dphisurj} has a solution is equivalent to show that 
\begin{equation}\label{eq:dphisurj2}
\int_0^T  U_0(t)^*(\delta H_0+\varepsilon(t)\delta\mu)U_0(t)dt=V,
\end{equation}
has a solution, with $V:=i U_0(T)^*V'\in\cal S$ since $V'\in \cal A_{H_0,\mu}$. 
A nice property of the trajectory $t\mapsto U_0(t)$ is that Equation~\eqref{eq:dphisurj2}
 can be solved explicitly. Indeed, let us denote by
$v_{a,b}$, $h_{a,b}$ and $m_{a,b}$, with $ a,b\in\N,\ 1\leq a,b\leq N,$
the coefficients of the matrices $V_0VV_0^*$, $V_0\delta H_0V_0^*$
and $V_0\delta\mu V_0^*$ respectively.
Expanding~\eqref{eq:dphisurj2} gives rise, in the case $a\neq b$ to 
\begin{eqnarray*}
v_{a,b}&=& h_{a,b}\int_0^T e^{i(\lambda_a-\lambda_b)(t-\frac
  T2) }dt\nonumber\\&&+m_{a,b}\int_0^T
\varepsilon(t)e^{i(\lambda_a-\lambda_b)(t-\frac T2)}dt\nonumber\\
&=&h_{a,b}\dfrac{\sin(\delta\lambda_{a,b}\dfrac T2)}{\delta\lambda_{a,b}}
+m_{a,b} \widehat\varepsilon (\delta\lambda_{a,b}),\label{eq:vab}
\end{eqnarray*} 
where $\delta\lambda_{a,b}=\lambda_a-\lambda_b$ and  $\widehat\varepsilon (\delta\lambda_{a,b})=\int_0^T
\varepsilon(t)e^{i\delta\lambda_{a,b}(t-\frac T2)}dt=\hat{\varepsilon}^r_{a,b}+i\hat{\varepsilon}^i_{a,b}$.\\In the case $a=b$, one finds that 
\begin{equation*}
v_{a,a}=h_{a,b}\dfrac T2
+m_{a,b} \widehat\varepsilon (0)=h_{a,b}\dfrac T2 +m_{a,b}
\int_0^T\varepsilon(t) dt.\label{eq:vaa}
\end{equation*}
 Note that the assumption $\delta H_0,\ \delta \mu\in \cal S_\R$
 combined with $\widehat\varepsilon
 (\delta\lambda_{a,b})=\overline{\widehat\varepsilon
   (\delta\lambda_{b,a})}$ implies that $v_{a,b}=\bar v_{b,a}$, so that $V\in\cal S$. The result follows.
\end{proof}
\begin{remark} 
In this theorem, we have defined $m_{a,a}$
arbitrarily. 
\end{remark}
This theorem gives a first hint about conditions required to identify
$(H_0,\mu)$. Condition~\eqref{cond1} is weaker to the standard
non-degeneracy condition $$\forall (a,b)\neq(a',b'),\ \lambda_b-\lambda_a\neq
\lambda_{b'}-\lambda_{a'},$$ 
and is in practice often satisfied. Condition~\eqref{cond2} deals
with the laser field itself. It is a non-resonant
condition to control the system.

\section{Numerical methods}\label{sec:algo}

In this section, we present two algorithms to
solve~\eqref{eq:target}. The strategy we follow is a direct adaptation
of previous results and proofs: we consider local approximations based on
 fixed point iterative solvers. In our approach, a crucial step
 consists in obtaining an appropriate time discretized version of~\eqref{eq:schrod1}.  
In the first part, we build such an approximation that enables the exact computation
of the derivative of the final state $U(T)$ with respect to $(H_0,\mu)$ and
derive from this setting a numerical strategy.

\subsection{Time discretization}

In order to simulate numerically Equation~\eqref{eq:schrod1}, we
introduce the following time discretization: give $N_T\in \N$, we denote
by $\Delta T=\dfrac T{N_T}$ the time step and for $n=0,\cdots,N_T$ by $U_n$ and $\varepsilon_n$
the approximations of $U(n\Delta T)$ and $\varepsilon(n\Delta T)$. In
order to preserve the unitary property of the matrices $U(t)$ at the
discrete level, we use a Crank-Nicholson scheme ruled by the formula:
$$i\dfrac{U_{n+1}-U_n}{\Delta T}=(H_0+\varepsilon_n\mu)\dfrac{U_{n+1}+U_n}2.$$
The corresponding iteration is then given by:
$$(Id+L_n)U_{n+1}=(Id-L_n)U_{n},
$$
where $L_n=\frac{i\Delta T}2 (H_0+\varepsilon_n\mu)$.\\
Let us now detail the effect of variations $\delta H_0$, $\delta \mu$
in $H_0$ and $\mu$ on the sequence $(U_n)_{n=0,...,N_T}$. We have:
\begin{eqnarray*}
(Id+L_n)\delta U_{n+1}+\delta L_nU_{n+1}&=&(Id-L_n)\delta U_{n}\\&&-\delta L_nU_{n},\\
\delta L_n(U_{n+1}+U_{n})&=&(Id-L_n)\delta U_{n}\\&&-(Id+L_n)\delta
U_{n+1},\\
(U_{n+1}+U_{n})^*\delta L_n(U_{n+1}+U_{n})&=&-2\Large(U_{n+1}^*\delta U_{n+1}\\&&-U_n^*\delta U_{n}\Large),
\end{eqnarray*}
where $\delta L_n=\frac{i\Delta T}2 (\delta H_0+\varepsilon_n\delta\mu)$.
This finally gives rise to:
\begin{align*}
U_{n+1}^*\delta U_{n+1}-U_n^*\delta U_{n}\phantom{\hspace{.6\linewidth}}
\\=-i\Delta T\dfrac{(U_{n+1}+U_{n})^*}2 (\delta H_0+\varepsilon_n\delta \mu)
\dfrac{U_{n+1}+U_{n}}2 .
\end{align*}
Since the initial value is fixed, we obtain:
\begin{eqnarray}
U_{N_T}^*\delta U_{N_T}\phantom{\hspace{.73\linewidth}}\nonumber\\
\!\!=\!\!-i\Delta T\!\!\sum_{n=0}^{N_T-1}\dfrac{(U_{n+1}+U_{n})^*}2 (\delta
H_0+\varepsilon_n\delta \mu)\dfrac{U_{n+1}+U_{n}}2.
\nonumber\\
\label{eq:ident}
\end{eqnarray}
This result can be seen as a discretized version of~\eqref{eq:idrem}
where $\mu$ is not necessarily null. We insist on the fact that such a result is specific to the
Crank-Nicholson discretization. As far as we know, no other numerical solvers
give rise to discretization of~\eqref{eq:idrem} where the
variations $\delta H_0$ and $\delta \mu$ are explicit.  
\subsection{Fixed points methods}\label{sec:fpm}
We now present some iterative solvers to compute solutions
of~\eqref{eq:target}. 

\subsubsection{A Newton Method}

In the
discrete setting, we still
denote by $\varphi$ the operator:
\begin{eqnarray*}
\varphi:& \cal S_\R\times \cal S^0_\R\rightarrow \cal U\\
&(H_0,\mu)\mapsto U_{N_T}.
\end{eqnarray*} 
To solve the equation $\varphi(H_0,\mu)=U_{target}$, a Newton method
would consist in the following iteration:
\begin{equation}\label{eq:newton}
d\varphi(H_0^k,\mu^k)\cdot(\delta
H_0^k,\delta\mu^k)=-\left(\varphi(H_0^k,\mu^k)-U_{target}\right),
\end{equation}
where $k$ is the iteration index, $\delta H_0^k=H_0^{k+1}-H_0^k$,
$\delta\mu^k=\mu^{k+1}-\mu^k$. \\
In our case, \eqref{eq:newton} reads:
\begin{equation*}
\delta U_{N_T}^k=U_{target}-U^k_{N_T}.
\end{equation*}
Using~\eqref{eq:ident}, one can rewrite this equation as follows:
\begin{eqnarray*}
\Delta T\sum_{n=0}^{N_T-1}\dfrac{(U^k_{n+1}+U^k_{n})^*}2 (\delta
H^k_0+\varepsilon_n\delta \mu^k)
\dfrac{U^k_{n+1}+U^k_{n}}2\\=i\left((U_{N_T}^k)^*U_{target}-Id\right), 
\end{eqnarray*}
where we recall that the unknowns are $\delta H^k_0$ and $\delta \mu^k$.
This equation has generally no solutions, since its left hand side
belongs to $\cal S$ what is not the case for its right hand side. To
solve this problem, we replace $i\left((U_{N_T}^k)^*U_{target}-Id\right)$ by a first
order approximation $S^k\in {\cal S_\R}$. Two possible choices are:
\begin{eqnarray}
\exp(-iS^k)&:=&(U_{N_T}^k)^*U_{target} \label{def:Ak1}\\
S^k&:=&i\dfrac{(U_{N_T}^k)^*U_{target}-U^*_{target}U_{N_T}^k}2.\label{def:Ak2} 
\end{eqnarray}
In the numerical
tests, the same behavior is observed when
using the first or the second definition.

\begin{remark}
The previous method can be simplified to obtain a procedure where no
matrix needs to be assembled and the 
inverted during iterations. Instead of up-dating at each iteration in the pair
$(H_0,\mu)$ in the term $d\varphi(H_0,\mu)$ of
Formula~\eqref{eq:newton}, one can  
keep a constant approximation $(H^{ref}_0,\mu^{ref})$ of the solution. We
denote by $(U^{ref}_n)_{n=0,\cdots,N_T}$ the corresponding sequence of states. The
iteration then reads:
\begin{eqnarray*}
\Delta T\sum_{n=0}^{N_T-1}\dfrac{(U^{ref}_{n+1}+U^{ref}_{n})^*}2 (\delta
H^k_0+\varepsilon_n\delta \mu^k)
\dfrac{U^{ref}_{n+1}+U^{ref}_{n}}2\\=S^k,
\end{eqnarray*}
where $S^k$ is defined in the previous section, see~\eqref{def:Ak1}
and~\eqref{def:Ak2}. Note that such an algorithm is actually a
time-discretized version of the fixed  
point used in the proof of Theorem~\ref{Th:ontoloc}, except that here
$\mu$ is not supposed to be null. 
\end{remark}

\subsubsection{Implementation of the iterative solvers}
Both previous methods require inversions of linear systems which are
not given explicitly in our formulations. To fill in this gap, we explain
here how to assemble the matrices, i.e. to rewrite the equation
\begin{equation*}
\Delta T\sum_{n=0}^{N_T-1}\dfrac{(U_{n+1}+U_{n})^*}2 (\delta
H_0+\varepsilon_n\delta \mu)
\dfrac{U_{n+1}+U_{n}}2=S,
\end{equation*}
in terms of linear system. In what follows,
we denote by $X_M$ the vector representation of a matrix 
$M$ consisting in concatenating vertically its columns. 
A first step to do this is to note that the later
equation reads as follows:
\begin{eqnarray}
\!\!\Delta T\!\left(\sum_{n=0}^{N_T-1}\!\! M_{U_{n+1/2}}\right)\! X_{\delta
  H_0}\!\!+\!\!\Delta T\!\left(\sum_{n=0}^{N_T-1}\!\! \varepsilon_n
M_{U_{n+1/2}}\right)\!X_{\delta \mu}\nonumber\\=X_S,\label{eq:Matrices}
\end{eqnarray}
with
$$M_{U_{n+1/2}}= kron({\bf 1}_{N_d},U_{n+1/2}^*).\times
kron(U_{n+1/2}^T,{\bf 1}_{N_d}).$$
Here, $kron$ denotes the Kronecker product, $U_{n+1/2}=\frac{U_{n+1}+U_{n}}2$ , the term by term 
product of two matrices $A$ and $B$ is denoted by $A .\times B$ and  ${\bf 1}_{N_d}$ denotes
the matrix of $\R^{N_d,N_d}$ whose coefficients are equal to
1.

A second step must then be carried out: since the matrices $\delta H_0$ and
$\delta \mu$ are symmetric, one has to consider the columns of the
matrices in \eqref{eq:Matrices} that correspond to the coefficients of $\delta H_0$ located,
e.g., above the diagonal and
the coefficients of $\delta \mu$ located strictly above the
diagonal. In the same way, only the lines of the resulting system that
correspond to the coefficients located above the diagonal of $S$
shall be considered. 

Taking the real and the imaginary part of the equations, the resulting
system is of size $N_d^4$.
\subsection{A continuation method for global controllability}\label{sec:glob}
The algorithms proposed in Section~\ref{sec:fpm} are only locally
convergent. The purpose of this section is to present a continuation
method that enables to extend their range of application.

As mentioned above, numerous methods exist to solve the control
problem where the laser term $\varepsilon$ in Equation~\eqref{eq:schrod1} is unknown and $H_0$ and $\mu$ are
given~\cite{MST,MT,BW}. Based on this fact, the method we propose is the following.
Given an initial guess $(H^0_0,\mu^0)$, find a control $\varepsilon^0$
such that $U^0_{N_T}$, the final state associated to $(H^0_0,\mu^0)$
correctly approximates $U_{target}$. Given $\theta\in[0,1]$ , we define the
interpolated fields
$\varepsilon^\theta=(1-\theta)\varepsilon^0+\theta\varepsilon$. A
fixed point method as the one presented in Section~\ref{sec:fpm} can then be
applied with  $(H^0_0,\mu^0)$ as an initial guess to solve the
operator control problem with $\varepsilon^\theta$. 
Our algorithm consists in repeating this
procedure by solving iteratively the operator control problem
associated to the field $\varepsilon^{k\delta\theta}$
using $(H^{k-1}_0,\mu^{k-1})$ as initial guess.
Carrying this procedure up to $\theta=1$ enables to solve the original problem.


\section{Numerical results}\label{sec:test}
In this last section, we present some numerical results obtained with
the algorithms of the previous sections. 
As a laser term in Equation~\eqref{eq:schrod1}, we use $\varepsilon(t)=\sin(t)$. The other numerical data
are $N_d=5$, $T_0=10$, $N_T=10^2$,  $T=2\pi T_0$ and $\Delta T=T/N$.


\subsection{Newton Method}
We first test our Newton method. In this way, we choose randomly a
pair $(H_0,\mu)$, with coefficients in $[-1,1]$ and compute the corresponding final state
$U_{N_T}$. Then, we start the Newton procedure with an initialization
$(H_0+\Delta H_0,\mu+\Delta \mu)$ where $(\Delta H_0,\Delta \mu)$ are
also chosen randomly. An example of computation is given in the next table.

\begin{table}[h!]
\begin{center}
\begin{tabular}{|c|c|c|}
\hline Iteration & $\log_{10}(\|H_0^k-H_0\|_{\cal U})$ & $\log_{10}(\|\mu^k-\mu\|_{\cal U})$\\
\hline 1&-1.579029&-1.358376\\
\hline 2&-3.003599&-2.865026\\
\hline 3&-4.339497&-4.122528\\
\hline 4&-8.234980&-8.179398\\
\hline 5&-13.963299&-14.029020\\
\hline 6&-14.022486&-14.131066\\
\hline 
\end{tabular}
\end{center}
\end{table}

Here, we refind a pair $(H_0,\mu)$ starting from a $10\%$ random
perturbation. We see that the numerical convergence is obtained after
6 iterations. Note also that the quadratic convergence is observed.  
\subsection{Continuation method}
In a second test, we use the continuation method presented
in~\ref{sec:glob} to tackle a problem where the algorithms of
Section~\ref{sec:fpm} do not apply. Given a target $U_{target}$ obtained with the
field $\varepsilon$ and a pair $(H_0,\mu)$ that is chosen randomly,
we look for the operators $H'_0$ and $\mu'$ that solve the
control problem associated to the field $\cos(3t)$ and the
target $U_{target}$.\\
The direct use of the Newton method of Section~\ref{sec:fpm} does not
work: in this case, the algorithm does not converge. The
continuation method enables to solve this problem. Using $\delta
\theta=1/4$, and 10 iterations of the Newton method as inner loop,
a relevant pair $(H_0',\mu')$ is obtained. \\
This example has been reproduced for numerous random initial pairs $(H_0,\mu)$.

\addtolength{\textheight}{-20cm}   








\bibliographystyle{IEEEtran}



\end{document}